\documentclass[12pt,thmsa]{article}
\usepackage{amsfonts}
\usepackage{sw20bams}

\input tcilatex
\begin{document}

\author{Steven Finch}
\title{A\ Convex Maximization Problem:\ Continuous Case}
\date{December 5, 1999}
\maketitle

\begin{abstract}
We study a specific convex maximization problem in the space of continuous
functions defined on a semi-infinite interval. An unexplained connection to
the discrete version of this problem is investigated.
\end{abstract}

\section{Problem}

Consider the class $X$ of continuous functions $f:[1,\infty )\rightarrow 
\Bbb{R}$ satisfying 
\[
\begin{array}{cccc}
f(x)\geq x &  & \text{for} & 1\leq x<2
\end{array}
\]
and

\[
\begin{array}{ccccc}
(y+1)f(y)+f(x)\geq (y+1)x &  & \text{for} & x\geq 2, & 1\leq y\leq x\text{ .}
\end{array}
\]
Prove that:

\begin{enumerate}
\item[(i)]  $x^{c\text{ }}\epsilon $ $X$ if and only if $c\geq \frac{1+\sqrt{%
5}}2$ (thus $X$ is nonempty)

\item[(ii)]  $\sup\limits_{f\text{ }\epsilon \text{ }X}\dint\nolimits_1^%
\infty \dfrac 1{f(x)}dx<\infty $

\item[(iii)]  the supremum in (ii) is attained by the function $a^{\text{ }%
}\epsilon $ $X$ defined by

\[
a(x)=\left\{ 
\begin{tabular}{llll}
$x$ &  & if & $1\leq x<2$ \\ 
$2\,(x-1)$ &  & if & $2\leq x<3$ \\ 
$3\,(x-2)$ &  & if & $5\leq x<8$%
\end{tabular}
\right. 
\]
\end{enumerate}

on $[1,3)\cup [5,8)$ and elsewhere by

\[
a(x)=(y+1)(x-a(y)) 
\]

where $y$ satisfies $(y+1)a^{\prime }(y)+a(y)=x$ and $a^{\prime }(y)$ is the
derivative of $a(y)$.\\

\noindent \textbf{Remark. }Part (iii) is, in essence, the continuous analog
of a certain number theoretic conjecture due to Levine and O'Sullivan \cite
{LOS} .

\section{Partial Solution}

After proving part (i), our treatment will be brief and rather informal. We
do not prove that the conjectured maximizing function $a(x)$ is well-defined
nor that it is feasible ($a$ $\epsilon $ $X$). Our purpose is to compute the
function $a(x)$ as far as possible (assuming it makes sense!), and to
describe a link between $a(x)$ and the discrete version $a_n$ studied in a
companion paper.

\subsection{Proof of (i)}

Define a family of functions 
\[
\Phi _y(x)=(y+1)\,y^c+x^c-(y+1)\,x 
\]
\noindent for all $x\geq 0$ and $y\geq 1$. It is sufficient to show that

\begin{enumerate}
\item[(a)]  given $1<c<\frac{1+\sqrt{5}}2$, there exist $x\geq 2$ and $1\leq
y\leq x$ with $\Phi _y(x)<0$

\item[(b)]  given $c=\frac{1+\sqrt{5}}2$, $\Phi _y(x)\geq 0$ for all $x\geq
0 $ and $y\geq 1$.
\end{enumerate}

\noindent Observe that, for fixed $y\geq 1$, $\Phi _y$ is minimized (via
calculus) at the point

\[
x_y=\left( \dfrac{y+1}c\right) ^{\frac 1{c-1}} 
\]
\noindent with minimum value 
\[
\Phi _y(x_y)=(y+1)\left[ y^c+x_y\,\left( \frac 1c-1\right) \right] <\Phi
_y(0) 
\]
\noindent (since $c>1$). Part (a) follows because 
\[
\lim\limits_{y\rightarrow \infty }\Phi _y(x_y)=-\infty 
\]
\noindent (since $1/(c-1)>c$). Part (b) follows by contradiction, because if 
$\Phi _y(x_y)<0$, then 
\[
y^c-\left( \dfrac{y+1}c\right) ^c\frac 1{c^2}<0 
\]
\noindent (since $1/(c-1)=c$), which implies 
\[
\left( \dfrac{y+1}y\right) ^c>c^{c+2} 
\]
\noindent that is, 
\[
y<\frac 1{c^{(c+2)\,(c-1)}-1}<1 
\]
\\\noindent But this is contrary to the hypothesis that $y\geq 1$. QED.

\subsection{Comment on (ii)}

A proof is not presently known, although it is expected to follow in a
manner similar to that in the discrete case \cite{Erdos}, \cite{LOS}.

\subsection{Comments on (iii)}

Again, a proof is not presently known. Calculus allows us, however, to
recursively unwrap the ''self-generating'' nature of the function $a(x)$ to
obtain some useful formulas. For example, if $1\leq y<2$, then $a(y)=y$, $%
a^{\prime }(y)=1$ and $(y+1)+y=x$; hence $x=2\,y+1$ and $y=$ $\frac
12\,(x-1) $. We deduce that $a(x)=\left( \frac 12\,(x-1)+1\right) \left(
x-\frac 12\,(x-1)\right) =\frac 14\,(x+1)^2$ for $3\leq x<5$. \\

Likewise, if $2\leq y<3$, then $a(y)=2\,(y-1)$, $a^{\prime }(y)=2$ and $%
2\,(y+1)+2\,(y-1)=x$; hence $x=4\,y$ and $y=$ $\frac 14\,x$. We deduce that $%
a(x)=\left( \frac 14\,x+1\right) \left( x-2\,(\frac 14\,x-1)\right) =\frac
18\,(x+4)^2$ for $8\leq x<12$. Proceeding similarly, the following is
obtained: 
\[
a(x)=\left\{ 
\begin{array}{cccc}
x &  & \text{if} & 1\leq x<2 \\ 
2\,(x-1) &  & \text{if} & 2\leq x<3 \\ 
\dfrac 14\,(x+1)^2 &  & \text{if} & 3\leq x<5 \\ 
3\,(x-2) &  & \text{if} & 5\leq x<8 \\ 
\dfrac 18\,(x+4)^2 &  & \text{if} & 8\leq x<12 \\ 
4\left( \dfrac x3\right) ^{3/2} &  & \text{if} & 12\leq x<27 \\ 
\dfrac 1{12}\,(x+9)^2 &  & \text{if} & 27\leq x<45 \\ 
\dfrac 14\left[ 1-8\,x+\left( \dfrac{8\,x+3}3\right) ^{3/2}\right] &  & 
\text{if} & 45\leq x<84 \\ 
(z^2+1)\left( x-\dfrac{4\,z^3}{3^{3/2}}\right) &  & \text{if} & 84\leq x<276
\\ 
\dfrac 4{81}\left[ 64-108\,x+\left( 9\,x+16\right) ^{3/2}\right] &  & \text{%
if} & 276\leq x<657 \\ 
\left( \dfrac{w^2+5}8\right) \left( x-1+\dfrac{w^2}4-\dfrac{w^3}{4\cdot
3^{3/2}}\right) &  & \text{if} & 657\leq x\leq 1781
\end{array}
\right. 
\]

\noindent where the auxiliary variables $z$ and $w$ are defined by 
\[
z=\dfrac 1{20^{1/3}}\left[ \left( \frac{135\,x^2+16}5\right)
^{1/2}+3^{3/2}x\right] ^{1/3}-\dfrac{20^{1/3}}5\left[ \left( \frac{%
135\,x^2+16}5\right) ^{1/2}+3^{3/2}x\right] ^{-1/3} 
\]

\noindent and 
\begin{eqnarray*}
w &=&\frac 15\,\left[
10\,(2700\,x^2+2106\,x+289)^{1/2}+3^{1/2}(300\,x+117)\right] ^{1/3} \\
&&+\frac{23}5\,\left[
10\,(2700\,x^2+2106\,x+289)^{1/2}+3^{1/2}(300\,x+117)\right] ^{-1/3}+\frac{%
4\cdot 3^{1/2}}5
\end{eqnarray*}

\noindent We have not attempted to determine $a(x)$\ for $x>1781$.

\subsection{Alternative Expression}

A more compact, but less explicit formula for $a(x)$ is as follows: 
\[
a(x)=\left\{ 
\begin{array}{cccc}
x &  & \text{if} & 1\leq x<2 \\ 
\max\limits_{1\leq y<x}\,(y+1)(x-a(y)) &  & \text{if} & x\geq 2
\end{array}
\right. 
\]
\noindent For example, if $2\leq x<3$, then 
\[
\max\limits_{1\leq y\leq 2}\,(y+1)(x-y)=2\,(x-1) 
\]
\noindent since the maximum cannot occur at $y=$ $\frac 12\,(x-1)<1$, hence
it must occur at one of the endpoints $y=1$ or $y=2$. Since $%
2\,(x-1)>3\,(x-2)$ for $x<4$, the claim is true. Suppose now that 
\[
a(x)=\max\limits_{1\leq y<x}\,(y+1)(x-a(y))>2\,(x-1) 
\]
\noindent This, in turn, implies that 
\[
2\,(x-1)<\max\limits_{2\leq y<x}\,(y+1)(x-a(y))\leq \max\limits_{2\leq y\leq
x}\,(y+1)(x-2\,(y-1))=3\,(x-2) 
\]
\noindent because the maximum cannot occur at $y=$ $\frac 14\,x<1$. But $%
3\,(x-2)<2\,(x-1)$, which yields a contradiction. Therefore $a(x)=2\,(x-1)$
for $2\leq x<3$.

Likewise, if $3\leq x<5$, then 
\[
\max\limits_{1\leq y\leq 2}\,(y+1)(x-y)=\frac 14\,(x+1)^2 
\]
\noindent since the maximum here occurs at $y=$ $\frac 12\,(x-1)$ and $1\leq
y<2$; and if $5\leq x<8$, then 
\[
\max\limits_{1\leq y\leq 2}\,(y+1)(x-y)=3\,(x-2) 
\]
\noindent since the maximum cannot occur at $y=$ $\frac 12\,(x-1)>2$.
Similar \textit{reductio ad absurdum} reasoning gives 
\[
a(x)=\left\{ 
\begin{array}{cccc}
\dfrac 14\,(x+1)^2 &  & \text{if} & 3\leq x<5 \\ 
3\,(x-2) &  & \text{if} & 5\leq x<8
\end{array}
\right. 
\]
\noindent as was to be proved.

An equivalence\ proof applicable for all $x\geq 8$ is not known. It will be
necessary to demonstrate that subinterval maximums always occur at interior
points, that is, at points where the derivative $a^{\prime }(x)$ vanishes.

\section{Link to Discrete Case}

In a companion paper, we studied the infinite sequence $a_1,a_2,...$,
defined by 
\[
\begin{array}{ccccc}
a_1=1, &  & a_2=2, &  & a_3=4
\end{array}
\]
and, when $i\geq 4$, 
\[
a_i=(j+1)(i-a_j) 
\]
where $j$ satisfies $(j+1)(a_j-a_{j-1})+a_{j-1}\leq i\leq
(j+2)(a_{j+1}-a_j)+ $ $a_j$.

An analogous alternative expression 
\[
a_i=\left\{ 
\begin{array}{cccc}
1 &  & \text{if} & i=1 \\ 
\max\limits_{1\leq j<i}\,(j+1)(i-a_j) &  & \text{if} & i\geq 2
\end{array}
\right. 
\]
\noindent applies here (although in this case a rigorous equivalence proof
is known). Such structural similarity leads us to expect a vague connection
between the sequence $a_i$ and the function $a(x)$, but the precise nature
of the link is difficult to anticipate. We empirically observe that the
non-analytic points $k$ of $a(x)$, that is, the subinterval endpoints in the
definition of $a(x)$, are evidently all integers. Further, the value of $a_k$
apparently coincides with $a(k)$ at all such points: 
\begin{eqnarray*}
a_1 &=&1=a(1) \\
a_2 &=&2=a(2) \\
a_3 &=&4=a(3) \\
a_5 &=&9=a(5) \\
a_8 &=&18=a(8) \\
a_{12} &=&32=a(12) \\
a_{27} &=&108=a(27) \\
a_{45} &=&243=a(45) \\
a_{84} &=&676=a(84) \\
a_{276} &=&4704=a(276) \\
a_{657} &=&19044=a(657) \\
a_{1781} &=&93925=a(1781) \\
a_{12460} &=&2148412 \\
a_{49312} &=&19916344 \\
a_{245395} &=&?
\end{eqnarray*}
\noindent It is possible that this pattern breaks down at some stage beyond
our computational means. We conjecture that this is not the case: that
instead $a_k$ and $a(k)$ are equal for infinitely many integers $k$. This
intriguing correspondence between the discrete and continuous versions is
presently without explanation.


\begin{thebibliography}{9}
\bibitem{Erdos}  Erd\"os, P., Remarks on number theory, III.\ Some problems
in additive number theory, \textit{Mat. Lapok} 13 (1962) 28-38.

\bibitem{Finch}  Finch, S. R., A convex maximization problem, \textit{J.
Global Optimization }2 (1992) 419; also MathSoft Inc., website URL
http://www.mathsoft.com/asolve/convex/convex.html, 1999.

\bibitem{LOS}  Levine, E. and O'Sullivan, J., An upper estimate for the
reciprocal sum of a sum-free sequence, \textit{Acta Arithmetica} 34 (1977)
9-24.\\

\begin{tabular}{l}
Steven Finch \\ 
MathSoft Inc., 101 Main Street \\ 
Cambridge, MA, USA 02142 \\ 
\textit{sfinch@mathsoft.com}
\end{tabular}
\end{thebibliography}
\end{document}